\documentclass[12pt,a4wide]{article}
\usepackage{amsmath, amsfonts, amssymb, amsthm, euscript, amscd, latexsym}
\usepackage[pdftex,colorlinks=true,
                       pdfstartview=fitv,
                       linkcolor=blue,
                       citecolor=blue,
                       urlcolor=blue,
           ]{hyperref}
\def\eee#1{ \begin{equation} #1 \end{equation} }
\def\aa#1{ \begin{align*} #1 \end{align*} }
\def\aaa#1{ \begin{align} #1 \end{align} }
\def\mm#1{ \begin{multline*} #1 \end{multline*} }
\def\mmm#1{ \begin{multline} #1 \end{multline} }

\newtheorem{thm}{\sc Theorem}

\newtheorem{cor}{\sc Corollary}

\newtheorem{pro}{\sc Proposition}

\newcommand{\sss}{\scriptscriptstyle}

\newcommand{\eps}{\varepsilon}
\newcommand{\pl}{\partial}
\newcommand{\gt}{\geqslant}
\newcommand{\lt}{\leqslant}
\newcommand{\sub}{\subset}

\newcommand{\dl}{\delta}
\newcommand{\al}{\alpha}
\newcommand{\gm}{\gamma}
 
 \newcommand{\Dl}{\Delta}
 \newcommand{\la}{\lambda}
 \newcommand{\sg}{\sigma}
\newcommand{\dd}{\diagdown}
\newcommand{\om}{\omega}
\newcommand{\mc}{\mathcal}

\newcommand{\td}{\tilde}

\newcommand{\e}{{\sss E}}

\newcommand{\x}{\times}
\newcommand{\te}{\theta}
\newcommand{\mto}{\mapsto}

\newcommand{\E}{\mathbb E}

\newcommand{\W}{\mathbb W}

\newcommand{\C}{{\rm C}}

\newcommand{\scal}{{\rm scal}}

\DeclareMathOperator{\ind}{\mathbb I}
\newcommand{\lap}{\Delta}
\newcommand{\nab}{\nabla}
\newcommand{\fdot}{\,\cdot\,}
\def\Rnu{{\mathbb R}}

\def\ffi{\varphi}

\def\intl{\int\limits}

\def\com#1{}

\long\def\symbolfootnote[#1]#2{\begingroup%
\def\thefootnote{\fnsymbol{footnote}}\footnote[#1]{#2}\endgroup}

\begin{document}

\author{Evelina Shamarova}

 \date{}

 \title{ 
 Chernoff's theorem for backward propagators
 and applications to diffusions on manifolds
 }

 \maketitle

 \vspace{-11mm}

 {\small
  \begin{center}
   Centro de matem\'atica da Universidade do Porto.\\
  { \hspace{2mm} E-mail: 
  \href{mailto:evelinas@fc.up.pt}{evelinas@fc.up.pt}}
 \end{center}
 }

 \vspace{3mm}

 \begin{abstract}
The classical Chernoff's theorem is a statement about discrete-time approximations of
semigroups, where the approximations are consturcted as products of time-dependent contraction
operators strongly differentiable at zero.
We generalize the version of Chernoff's theorem for semigroups proved in \cite{sm_vr} (see also 
\cite{chernoff1} and \cite{chernoff}), and obtain
a theorem about discrete-time approximations of backward propagators. 
 \end{abstract}

\vspace{1mm}


\vspace{4mm}

  \noindent\textit{2010 Mathematics Classification:} {47D06, 47D07}

  %
  %

  \section{Introduction}
Let $E$ be a Banach space, and let $\mc L(E)$ denote the space of all bounded
operators $E\to E$. 
  Let $U(s,t)$ be a backward propagator on $E$ possessing the left generator $A_t$.
For convenience we give definitions of backward propagators and their left generators
(see \cite{bp}).
A two-parameter family of operators $\{U(s,t)\in \mc L(E)\, : \, 0\lt S\lt s \lt t \lt T\}$
is called a backward propagator on $E$ if
\aaa{
\label{flow}
&U(s,t) = U(s,\tau)U(\tau,t),\\
&U(s,s) = I\notag
} 
for all $s$, $\tau$, $t$ such that $S\lt s\lt \tau\lt t \lt T$.
The operator $A_t$ on $E$ defined as
\aa{
A_tx = \lim_{h\downarrow 0} \frac{U(t-h,t)x - x}{h},
}
$t>0$, with the domain $D(A_t)$ consisting of those $x\in E$ for which the above limit exists,
is called the left generator of the backward propagator $U(s,t)$.

Analogously, we can introduce the concept of the right generator of 
a backward propagator (see  \cite{bp}):
the operator $A^+_t$ on $E$ defined as
\aa{
A^+_tx = \lim_{h\downarrow 0} \frac{U(t,t+h)x - x}{h},
}
$t>0$, with the domain $D(A^+_t)$ consisting of those $x\in E$ for which the above limit exists,
is called the right generator of the backward propagator $U(s,t)$.

   Let $Q_{s,t}$, $0\lt S\lt s \lt t \lt T$, be a two-parameter family of contraction operators on  $E$,
   whose left derivatives at $s=t$ equal to $A_t$. 
   The discrete-time approximations of the backward propagator $U(s,t)$ are constructed as
   products of $Q_{t_1,t_2}$, $s\lt t_1 \lt t_2 \lt t$.
Note that we could equivalently use right generators
of the backward propagator and right derivatives of $Q_{s,t}$ at $t=s$.
The theorem will also work in the situation
with (forward) propagators and
the two-parameter family of contractions $Q_{t,s}$ parametrized by times
$t$ and $s$ such that $T \gt t \gt s \gt S \gt 0$.
We prove our main result for backward propagators because
in the application to diffusions on manifolds (Section \ref{3.2}) 
backward propagators will be associated to transition probability functions.
Specifically, we consider the situation when the backward propagator is represented
   by a transition probability function of a time-inhomogeneous diffusion on a compact Riemannian manifold,
   the contraction operators are integral operators with probabilistic kernels,
   the left generators of the backward propagator are second-order differential operators on the manifold,
   and the discrete-time approximations are distributions of diffusion processes in the surrounding
   Euclidean space. 
     We then obtain the approximation of the distribution on the manifold by
   distributions in the Euclidean space.
   
Compared to the situation considered in \cite{T}, 
the stochastic processes under consideration are non-homogeneous. In particular,
the coefficients of the second-order differential operator representing the generator of 
the manifold-valued diffusion are time-dependent. Therefore, in the current paper we consider a more
general situation compared to \cite{sm_vr}, \cite{chernoff1}, \cite{chernoff}, and \cite{T} for both Chernoff's theorem
and its applications to diffusions on manifolds.

   \renewcommand{\labelenumi}{\theenumi)}

\section{Chernoff's theorem for backward propagators}
 \begin{thm}[Chernoff's theorem for backward propagators]
 \label{chernnon}
   Let  $U(s,t)$, $0\lt S\lt s \lt t \lt T$, be a backward propagator
   with the left generators $A_t$, and let $Q_{t_1,t_2}$, $S \lt t_1 \lt t_2 \lt T$, 
   be a two-parameter family of contractions $E\to E$.
   We assume that the following assumptions are fulfilled:
\begin{enumerate}
\item
\label{domain}
The subset $\cap_{t\in [S,T]}D(A_t)$ is dense in $E$.
\item
\label{norm}
 There exists a dense in $E$ Banach space $Y$ such that
 $Y\sub \cap_{t\in [S,T]}D(A_t)$ and
 $U(s,t) Y\sub Y$ for all $s,t \in [S,T]$, $s<t$, and, moreover, so that 
 there exists a constant $\gm > 0$ such that the norm in $Y$
satisfies the inequality $\|x\|_Y \gt \gm [\,\|x\|_\e + \sup_{\tau\in [S,T]}\|A_\tau x\|_\e\,]$.
\item
\label{continuous}
 For every $x\in Y$ and $t\in [S,T]$,
 the function $[S,t] \to Y$, $\tau\mto U(\tau,t)x$ is continuous.
\item
\label{continuous1}
For every $x\in Y$, the function $[S,T]\to E$, $t\mto A_tx$ is continuous.
\item
\label{limit}
 For all $x\in Y$ there exists the uniform in $t$ limit
 \aa{
 \lim_{h\downarrow 0} \frac{Q_{t-h,t}x - x}{h} = A_t x.
 }
\end{enumerate}
 Then,
 for any subinterval $[s,t]\sub [S,T]$,
 for any sequence of partitions 
 $\{s=t_0 < t_1 < \cdots < t_n=t\}$
such that $\max \,|t_{j}-t_{j-1}|\to 0$ as $n\to\infty$,
 and for all $x\in E$,
\aa{
Q_{t_0,t_1}\ldots Q_{t_{n-1},t_n}x
 \to U(s,t)\,x, \quad n\to \infty.
 }
 \end{thm}
 \begin{proof}
 Fix an $x\in Y$. First consider the case $s>S$.
 Using relation \eqref{flow}, we obtain:
 \mmm{
 \label{4}
 Q_{t_0,t_1}Q_{t_1,t_2}\ldots Q_{t_{n-1},t_n}-
 U(s,t) \\ =
 \sum_{j=1}^{n}Q_{t_0,t_1}\ldots Q_{t_{j-2},t_{j-1}}
 (Q_{t_{j-1},t_{j}}-U(t_{j-1},t_j))
 U(t_j,t).
 }
 Let $\dl_n=\max_j(t_j-t_{j-1})$, $j\gt 1$,
 be the mesh of the partition $\{s = t_0 < t_1 < \cdots < t_n=t\}$.
 Relation \eqref{4} implies:
 \aaa{
 \notag
 &\quad\|(Q_{t_0,t_1}Q_{t_1,t_2}\ldots Q_{t_{n-1},t_n}-
 U(s,t))x\|_\e\\
 \notag
 &\lt
 \sum_{j=1}^{n}\Dl t_j
 \Bigl\|
 \Bigl(
 \frac{Q_{t_{j-1},t_j}-I}{t_j-t_{j-1}}
 -\frac{U(t_{j-1},t_j)-I}{t_j-t_{j-1}}
 \Bigr)U(t_j,t)\,x
 \Bigr\|_{\sss E}\\
 \notag
 &\lt
 (t-s)\sup
 \Bigl\{
 \Bigl\|
 \Bigl(
 \frac{Q_{\tau-h,\tau}-I}{h}
 -\frac{U(\tau-h,\tau)-I}{h}
 \Bigr)
 U(\tau,t)\, x
 \Bigr\|_{\sss E}:\\
 \notag
 &\phantom{(t-s)\sup
 \Bigl\{
 \Bigl\|
 \Bigl(
 \frac{Q_{\tau-h,\tau}-I}{h}
 -\frac{U(\tau-h,\tau)-I}{h}
 \hspace{8mm}
 }
 \tau \in (s,t], h\in (0,\dl_n) \Bigr\}\\
 &\lt
 \label{1term}
 (t-s)\sup
 \Bigl\{
 \Bigl\|
 \Bigl(
 \frac{Q_{\tau-h,\tau}-I}{h}
 -A_\tau
 \Bigr)
  U(\tau,t)\,x
 \Bigr\|_{\sss E}:
 \tau \in (s,t], h\in (0,\dl_n)
 \Bigr\}\\
 &+(t-s)\sup
 \Bigl\{
 \Bigl\|
 \Bigl(
 \frac{U(\tau-h,\tau)- I}{h}-A_\tau
 \Bigr)
 U(\tau,t)\,x
 \Bigr\|_{\sss E}:
 \label{3term}
 \tau \in (s,t], h\in (0,\dl_n)
 \Bigr\}.
 }
 Note that for every $x\in Y$,
\aaa{
\label{11}
\Bigl(\frac{U(\tau-h,\tau)- I}{h}-A_\tau\Bigr)x
}
 converges to zero uniformly in $\tau \in [s,t]$.
 Indeed, by Assumption \ref{continuous1}, one can find a $\theta \in (0,1)$ such that
 $\frac{U(\tau-h,\tau)\,x- x}{h} = A_{\tau -\theta h}x$.
 Since the function $[s,t]\to E$, $\tau \mto A_\tau x$
 is continuous by assumption, it is also uniformly
 continuous which implies the uniform convergence in \eqref{11}.
Let $B_{\tau-h,\tau}$ denote one of the operators
 $\frac{Q_{\tau-h,\tau}-I}{h}-A_\tau$ or $\frac{U(\tau-h,\tau)-I}{h}-A_\tau$.
 We know that for every $x\in Y$, $B_{\tau-h,\tau}x$ converges
 to zero uniformly in $\tau\in [s,t]$. We would like to prove
 that $B_{\tau-h,\tau}U(\tau,t)\,x$ also converges to zero
 uniformly in $\tau\in [s,t]$.
 By the continuity of the map $[s,t]\to Y, \; \tau\mto U(\tau,t)y$,
 the set $\{U(\tau,t)x, \; \tau\in [s,t]\}$
 is a compact in $Y$. 
 We fix an arbitrary small $\eps>0$ and
 find a finite $\eps$-net $\{y_i\}_{i=1}^N\sub Y$ for
 this compact.
 Let us consider now $B_{\tau-h,\tau}$ as an operator from $Y$ to the Banach
 space $\mc E$ of continuous functions $[s,t]\to E$ with the norm $\sup_{\tau\in [s,t]}\|f_\tau\|_\e$.
 By the Banach-Steinhaus theorem the norms $\|B_{\tau-h,\tau}\|_{\mc L(Y,\mc E)}$
 are bounded. This implies the uniform  in $\tau\in [s,t]$ convergence to zero of $B_{\tau-h,\tau}U(\tau,t)\,x$,
 and therefore the convergence to zero of terms \eqref{1term} and \eqref{3term}.

 Thus, we have proved that $Q_{t_0,t_1}\ldots Q_{t_{n-1},t_n}x
 \to U(s,t)x$ as $n\to \infty$ for each $x\in Y$
 where $Y$ is dense in $E$.
 Since the operators $Q_{t_0,t_1}\ldots
 Q_{t_{n-1},t_n}$ are contractions, the convergence $Q_{t_0,t_1}\ldots Q_{t_{n-1},t_n}x
 \to U(s,t)x$ holds for all $x\in E$ by the Banach-Steinhaus theorem.
 We proved the theorem for the case $s>S$.

 Let us consider the case $s=S$. Fix an $x\in Y$. Let $s_N>s$ be a decreasing sequence of real numbers
 such that $\lim_{N\to\infty} s_N = s$. 
 Consider a partition $\mc P_N = \{s_N< t_1 < \cdots < t_n = t\}$ of $[s_N,t]$. 
 We have:
   \mmm{
   \label{12}
  \|Q_{s,s_N}Q_{s_N,t_1} \ldots Q_{t_{n-1},t_n}x - U(s,t)x\|_ {\mc L(E)}\lt \\
  \|Q_{s,s_N}(Q_{s_N,t_1} \ldots Q_{t_{n-1},t_n}x - U(s_N,t)x)\|_ {\mc L(E)}
   + \|(Q_{s,s_N}-U(s,s_N))U(s_N,t)x\|_ {\mc L(E)}.
  }
Let us prove that as $N\to \infty$, $(Q_{s,s_N}-U(s,s_N))U(s_N,t)x\to 0$.
We have:
\mmm{
\label{890}
(Q_{s,s_N}-U(s,s_N))U(s_N,t)x  =   \bigl(Q_{s_N-(s_N-s),s_N}-I\bigr) U(s_N,t)\, x \\ - 
  \bigl(U(s_N-(s_N-s),s_N)-I\bigr)\, U(s_N,t)\, x.
}    
We have proved that for every $x\in Y$, $B_{\tau-h,\tau}U(\tau,t)x$ converges
to zero uniformly in $\tau\in [s,t]$. This implies that the both summands
in \eqref{890} converge to zero as $N\to\infty$.
Further note that as the mesh of $\mc P_N$ tends to zero,
 \eee{
  \label{08}
   Q_{s_N,t_1}\ldots Q_{t_{n-1},t_n} x - U(s_N,t)\,x \to 0,
  }
since we can repeat the argument that leads to estimates
\eqref{1term} and \eqref{3term}.
To make our argument precise, we define $U(\tau-h,\tau) = U(s,\tau)$ and 
  $Q_{\tau-h,\tau} = Q(s,\tau)$ if $\tau-h<s$.
Next, since $Q_{s,s_N}$ is a contraction, we conclude 
that the first summand in \eqref{12} converges to zero
as the mesh $|\mc P_N|$ goes to zero. 
Thus for any $x\in Y$, the left-hand side of \eqref{12} converges to zero.
By the Banach-Steinhaus theorem it converges to zero for all $x\in E$.
The theorem is proved. 
  \end{proof}
 \begin{cor}[The case of commuting generators]
 \label{cor}
 Let $A_t$ be a stable (see \cite{sch}) family of pairwise commuting generators of
 strongly continuous semigroups, and let
 $Q_{t_1,t_2}$, $t_1,t_2>0$, be a
 two-parameter family of contraction operators $E\to E$,
 such that Assumptions \ref{domain}--\ref{limit} of Theorem \ref{chernnon}
 are fulfilled.
 Then, for any subinterval $[s,t]\sub [S,T]$, for any sequence of
 partitions
 $\{s=t_0 < t_1 < \cdots < t_n=t\}$
 of $[s,t]$ such that
 $\max \,(t_{j+1}-t_j)\to 0$ as $n\to\infty$,
 and for all $x\in E$,
 \[
 Q_{t_0,t_1}\ldots Q_{t_{n-1},t_n}x
 \to e^{\int_s^t A_r dr}\,x, \quad n\to \infty.
 \]
 \end{cor}
 For the proof of Corollary \ref{cor} we will need Proposition \ref{pro1}
 below (see \cite{sch}, p.489 for details).
  \begin{pro}
 \label{pro1}
 Let $\{A_t\}$ be a stable family of pairwise commuting generators
 of strongly continuous semigroups. Let us assume that there exists a
 space $Y\sub \cap_{t\in [S,T]}D(A_t)$
 which is dense in $E$, and let for all $y\in Y$,
 the mapping $[S,T]\to E, \, t\mto A_ty$ be continuous. Then,
 $(\int_s^t A_r dr, Y)$
 is closable and its closure (which we
 still denote by $\int_s^t A_r dr$) is a generator. Moreover,
 the backward propagator with the left generator $A_t$
 takes the form:
 \[
 U(s,t)= e^{\int_s^t A_r dr}.
 \]
 \end{pro}
  \begin{proof}[Proof of Corollary \ref{cor}]
Proposition \ref{pro1} and Theorem \ref{chernnon} imply 
Corollary \ref{cor}.
 \end{proof}

 \section{Application to diffusions on manifolds}
 Let $M$ be a $d$-dimensional compact Riemannian manifold isometrically embedded into a Euclidean space $\Rnu^m$. 
 Further let $\sg(t)$ be a nongenerate matrix in $\Rnu^m$.
We assume that the map $[S,T] \to GL(m)$, $t\mto \sg(t)$ is continuously
differentiable, 
where $GL(m)$ denotes the space or real nongenerate matrices $m\x m$.
 Consider the transition density function
 \aaa{
 \label{df}
p(s,x,t,y) = \frac{\det{\sg(t)}}{(2\pi(t-s))^{\frac{m}2}}
\exp\Bigl( -\frac{\bigl|\sg(t)y-\sg(s)x\bigr|^2_{\Rnu^m}}{2(t-s)}
\Bigr).
 }
One can easily verify that the non-homogeneous Markov process
associated to \eqref{df} is $x + \sg(t)^{-1}W_t$, where $W_t$ is
an $\Rnu^m$-valued Brownian motion.
\subsection{A short time asymptotic of a Gaussian-type integral operator}
In this paragraph we obtain a short time asymptotic for the intergral of the form
$\frac1{(2\pi t)^{\frac{d}2}}\int_M g(z)e^{-\frac{|z-y|^2}{2t}}
 \la_M(dz)$, where $\la_M$ is the volume measure on $M$.
 Unlike the short time asymptotic of the same integral obtained in \cite{sm_vr}
 we compute the coefficient at $t$ precisely. 
 In \cite{sm_vr}, the authors do not obtain the precise expression for this coefficient.  

 Let $\scal_M$ denote the scalar curvature, and  $\lap_M$ denote the Laplace-Beltrami operator on $M$.
 \begin{pro}
 \label{pr1}
 Let  $g\in C^2(M)$.
 Then, there exist a constant $K$, a time $t_0$, and a function $R: [0,t_0] \x M \to \Rnu$ satisfying
 $|R(t,y)|<Kt^{\frac12}$ for all $y\in M$ and for all $t\in [0,t_0]$ such that
 \mmm{
 \label{gauss_as}
 \frac1{(2\pi t)^{\frac{d}2}}\int_M g(z)e^{-\frac{|z-y|^2}{2t}}
 \la_M(dz)
 =g(y) -\frac{t}2\lap_M g(y) \\ 
 - g(y)\, \Bigl(\frac16\,\scal(y)+\frac1{16}\, \lap_M\lap_M\left.|\fdot - \, y|^2\right|_y\Bigr)\, t
  +t R(t,y)
 }
for all $y\in M$ and for all $t\in [0,t_0]$.
 \end{pro}
\begin{proof}
 Let $\imath$ be the isometrical embedding of $M$ into $\Rnu^m$.
 It is well known that
 $|\imath(z)-\imath(y)|^2=d(y,z)^2+\ffi(y,z)$,
 where $d$ is a geodesic distance in $M$, and
 $\ffi(y,z)=O(d(y,z)^4)$. 
 Let $U_y\sub M$ be a neighborhood of $y$, $U$ be a neighborhood
 of zero in the tangent space $T_y$ at $y$.
 Let $\psi_y: U\to U_y$
 be the diffeomorphism providing the normal
 coordinates in $U_y$, $f_y(x)=\ffi(y,\psi_y(x))$,
 $h_y(x)=\sqrt{\det g_{ij}(x)}\,g(\psi_y(x))$ where
 $g_{ij}$ is the metric tensor. We have:
 \[
 \int_{U_y} \hspace{-2mm} e^{-\frac{|z-y|^2}{2t}}g(z) \la_M(dz)
 =\int_{U_y}\hspace{-2mm} e^{-\frac{d(y,z)^2+\ffi(y,z)}{2t}}g(z) \la_M(dz)
 =\int_U e^{-\frac{|x|^2+f_y(x)}{2t}}h_y(x)dx.
 \]
 By results of~\cite{sm_vr}, there exist a function $\td R(t,\fdot)$ and a constant $\td K$ such that
 \eee{
 \label{asss}
 \frac1{(2\pi t)^{\frac{d}2}}
 \int_U e^{-\frac{|x|^2+f_y(x)}{2t}}h_y(x)dx
 =h_y(0)+\frac{t}2\,\lap h_y(0)-
 \frac{t}{16}\, h_y(0)\lap\lap f_y(0)+
 t\,\td R(t,x),
 }
 and $|\td R(t, \fdot)|<\td K t^{1/2}$.
 By arguments of \cite{sm_vr}, the neighborhood $U\sub \Rnu^d$ and the constant $\td K$ 
 can be chosen the same for all $y\in M$.
 Note that $h_y(0)=g(y)$. Next, it was obtained in~\cite{sm_vr} that
 $\lap h_y(0)=-\lap_M u(y)-\frac13\, u(y) \,\scal(y)$.
 Let us compute $\lap \lap f_y(0)$. Note that $\lap \lap\,
 d(y,\psi_y(x))^2=\lap \lap |x|^2=0$.
Hence, 
 \aa{
 \lap\lap f_y(0)=\lap\lap\bigl(|\imath \circ
 \psi_y(x)-\imath(y)|^2\bigr)|_{x=0}=\lap_M\lap_M \left.|\fdot - \, y|^2\right|_y.
 }
 Substitute  the expressions for $\lap h_y(0)$
 and $\lap\lap f_y(0)$ into~\eqref{asss}.
Next, we need to estimate the integral $\frac1{(2\pi t)^{\frac{d}2}}
 \int_{M\dd U_y} g(z) e^{-\frac{|z-y|^2}{2t}} \la_M(dz)$.
 Neighborhoods $U_y$ can be choosen  of the form
 $U_y=\{z\in M : |z-y|<\eps_y\}$ where $\eps_y$ can be choosen 
 bounded away from zero (see \cite{sm_vr}), say, by $\eps$.
 Let $t_0 > 0$ be small enough so that
 \eee{
 \label{linaaa}
 \frac1{(2\pi t)^{\frac{d}2}}
 \int_{M\dd U_y} g(z) e^{-\frac{|z-y|^2}{2t}} \la_M(dz)
 \lt \frac1{(2\pi t)^{\frac{d}2}} \, e^{-\frac{\eps^2}{2t}}
 \int_M |g(z)|\la_M(dz) < t^{3/2}
 }
 for $t < t_0$. 
  Estimate \eqref{linaaa} and the choice of the function $\td R$ imply 
 \eqref{gauss_as} with $R(t,y)$
 satisfying $|R(t,y)|<K t^{1/2}$, where the constant
 $K$ does not depend on $y$.
 \end{proof}
 \begin{cor}
 \label{cor1}
 Let $g\in C^2(M)$. Then,
 there exist a constant $K$, a time $t_0$, and a function $\bar R: [0,t_0] \x M \to \Rnu$ satisfying
 $|\bar R(t,x)|< Kt^{\frac12}$ for all $x\in M$ and for all $t\in [0,t_0]$ such that
for all $x\in M$, and for all $t\in [0,t_0]$,
 \aa{
 \frac{\int_M g(y)e^{-\frac{|y-x|^2}{2t}}\la_M(dy)}{\int_M e^{-\frac{|y-x|^2}{2t}}\la_M(dy)}  
 = g(x) -\frac{t}2\lap_M g(x) + t \bar R(t,x).
 }
 \end{cor}
 \begin{proof}
 The statement of the corollary easily follows from 
 Proposition \ref{pro1} applied to the functions $g(y)$ and
 $g(y)\equiv 1$ respectively.
 \end{proof}
\subsection{Surface measure generated by a non-homogeneous diffusion}
\label{3.2}
Consider the integral operator $\C(M)\to \C(M)$:
\aaa{
\label{Q}
(Q_{\tau-h,\tau}f)(x) = \frac{\int_M p(\tau-h,x,\tau,y) f(y)\la_M(dy)}{\int_M p(\tau-h,x,\tau,y)\la_M(dy)}.
}
After introducing the notation 
\aa{
p^M(\tau-h,x,\tau,y) =\frac{\ind_M(y) \, p(\tau-h,x,\tau,y)}{\int_{M}p(\tau-h,x,\tau,y) \la_M(dy)}
}
we can write \eqref{Q} in the form:
\aaa{
\label{QpM}
(Q_{\tau-h,\tau}f)(x) = \int_M p^M(\tau-h,x,\tau,y) f(y) \la_M(dy).
} 
Consider the operator product:
 \mmm{
 \label{Qprod}
(Q_{t_0,t_1}Q_{t_1,t_2}\ldots Q_{t_{n-1},t_n}f)(x)\\
= \int_M p^M(t_0,x,t_1,x_1) \la_M(dx_1)\int_Mp^M(t_1,x_1,t_2,x_2)\la_M(dx_2)\\
\ldots\int_M p^M(t_{n-1},x_{n-1},t_n,x_n)f(x_n)\la_M(dx_n).
 }
\begin{thm}
\label{thm_app}
Let the operator $Q_{\tau - h,\tau}: \C(M)\to \C(M)$ be defined by \eqref{QpM}. Then,
as the mesh of $\mc P$ tends to zero, the operator product defined by \eqref{Qprod}
converges at every point $f\in \C(M)$ with respect to the norm of $\C(M)$
to the backward propagator $U(s,t)$ whose left generator
is given by 
\aaa{
\label{generator}
(A_tf)(x) = -\frac12\, \lap_{M_t}f_t(\sg(t)x) + (\nab_{M_t}f_t(\sg(t)x),\sg'(t)x)_{\Rnu^m}, 
}
where $M_t = \sg(t) M$, $f_t = f\circ \sg(t)^{-1}$, $x\in M$. 
\end{thm}
\begin{proof}
Let us first show that the operators $A_t$ generate a non-homogeneous diffusion on $M$.
Note that $M_t$ is also isometrically embedded into $\Rnu^m$. The isometric
embedding $\imath_t$ defines a metric tensor $\td g_{ij}(t,x) = \sum_\al \frac{\pl\imath_t^\al}{\pl x^i}
\frac{\pl \imath_t^\al}{\pl x^j}(x)$ on $M_t$, and the Levi-Civita connection $\td \Gamma^i_{jk}(t,\fdot)$ of
the metric $\td g_{ij}(t,\fdot)$. Let $f\in \C^2(M)$ and 
$\td x = \sg(t) x\in M_t$. 
Further let $\{\td x_i\}$ be local coordinates in a neighborhood  $U$ of  $\td x$. 
We have:
\aa{
(A_tf)(x) = \td g^{ij}(t,\td x)\frac{\pl^2 f_t}{\pl \td x^i \pl \td x^j}(\td x)
- \td g^{ij}(t,\td x)\td\Gamma^k_{ij}(t,\td x)\frac{\pl f_t}{\pl \td x^k}(\td x)
+ \frac{\pl f_t}{\pl \td x^i}\, \widetilde{(\sg'(t)x)}^i
}
where $\widetilde{(\sg'(t)x)}^i$ are the coordinates, 
with respect to the basis $\frac{\pl}{\pl \td x^i}$, of the projection of the vector $\sg'(t)x$
onto the tangent space $T_{\td x}(M_t)$.
The matrix $\sg(t)^{-1}$ can be regarded as a change of coordinates in $U$. 
Let $\{x_i\} = \sg(t)^{-1}\{\td x_i\}$ be the new coordinates
in $U$ due to this change. Further let 
$g^{ij}(t,\fdot)$ and $\Gamma^k_{ij}(t,\fdot)$ denote the metric tensor and 
the Levi-Civita connection written in the coordinates $\{x_i\}$.
Note that $\{x_i\}$ are also local coordinates in the neighborhood $\sg(t)^{-1}U \sub M$
of the point $x \in M$. Also, $g^{ij}$ and $\Gamma^k_{ij}$
can be ragarded as the metric tensor and the Levi-Civita connection on $M$.
Taking into account this, we obtain 
the following connection between $g^{ij}$ and $\td g^{ij}$, $\Gamma^k_{ij}$ and $\td \Gamma^k_{ij}$:
\aa{
&\td g^{ij}(t,\td x) = g^{pq}(t, x) \sg^i_p(t) \sg^j_q(t),\\
&\td \Gamma^k_{ij}(t,\td x) = \sg_l^k(t) (\sg^{-1})_i^p(t) (\sg^{-1})_j^q(t) \Gamma^l_{pq}(t, x).
}
Moreover, $\frac{\pl^2 f_t}{\pl \td x^i\pl \td x^j}(\td x) =
 \frac{\pl^2 f}{\pl x^k\pl x^l}(x)
(\sg^{-1})^k_i(\sg^{-1})^l_j$ and $\frac{\pl f_t}{\pl \td x^k}(\td x)= 
\frac{\pl f}{\pl x^m}(x)(\sg^{-1})^m_k$.
This implies that
\aaa{
\label{lap2}
(A_t f)(x) = 
 g^{pq}(t,x)\frac{\pl^2 f}{\pl x^p\pl x^q}(x)- g^{pq}(t,x)\Gamma^k_{pq}(t, x)
\frac{\pl f}{\pl  x^k}(x) + \frac{\pl f}{\pl x^p}\,(\sg'(t)x)^p
}
where $(\sg'(t)x)^p$ are the coordinates, 
with respect to the basis $\frac{\pl}{\pl x^p}$, 
of the projection of the vector $\sg'(t)x$ onto the tangent space $T_x(M)$.
The existence of a unique diffusion on $\bigl(M, g_{ij}(t,\fdot)\bigr)$ generated
by the time-infomogeneous differential operator 
on the right-hand side of \eqref{lap2} is known (see, for example, \cite{takeyama}).
Therefore the operator $A_t$ defined by \eqref{generator} 
generates a diffusion on $M$. 

Let us show that Assumptions  \ref{domain}--\ref{limit} of Theorem \ref{chernnon} are fulfilled.
Note that all the generators $A_t$ have the same domain $\C^2(M)$. Therefore,
the space $Y$ can be taken to be the common domain $\C^2(M)$.
The norm in $Y$ is the following: $\|x\|_{\sss Y} = \|x\|_{\C(M)} + \sup_{\tau \in [S,T]} \|A_\tau x\|_{\C(M)}$.
Let $f\in \C^2(M)$, and let
$u(s,x)$ be the solution to the following final value problem on $\C(M)$:
\aaa{
\label{final2}
\begin{cases}
\frac{\pl u}{\pl s}(s,x) = -A_s \, u(s,x)\\
\lim_{s\uparrow t} u(s,x) = f(x).
\end{cases}
}
Further let $P(s,x,t,A)$ be the transition probability function of the diffusion
generated by $A_s$. Then, the backward propagator $U(s,t)$ can be 
expressed via $P(s,x,t,A)$:
\aa{
(U(s,t)f)(x) = 
\begin{cases}
\int_M P(s,x,t,dy) \, f(y), \; s<t, \\
f(x), \; s=t.
\end{cases}
}
Moreover, $u(s,x) = (U(s,t)f)(x)$ (see \cite{bp}).
Clearly, $u(s,\fdot)\in \C^2(M)$, and therefore Assumption  \ref{norm} is fulfilled.
Next, it is known that $u \in C^{2,1}(M\x[S,t])$, which implies that the map
$[S,t] \to C^2(M)$, $s\mto u(s,\fdot)$ is continuous. Therefore,
Assumption \ref{continuous} of Theorem \ref{chernnon} is also fulfilled. 

Let us show now that Assumption \ref{limit} is fulfilled.
Let $y_t = \sg(t)y$, $x_s = \sg(s)x$. Then 
\aa{
p(s,x,t,y) = \det \sg(t)\,q(t-s, x_s, y_t) 
}
where $q(\tau,x,y)$ is the Gaussian density  with respect to the Lebesgue measure on $\Rnu^m$. 
%
It is easy to verify that
\aaa{
\label{int.dens}
\int_M p(s,x,t,y) f(y) \, \la_M(dy) = \det \sg(t) \int_{M_t} q(t-s, x_s, y_t) \, 
f_t(y_t)\, \la_{M_t}(dy_t).
}
Using this formula and 
canceling the multiplier $e^{\frac{-|x_t-x_{t-\dl}|^2_{\Rnu^m}}{2\dl}}$ 
in the numerator and the denominator of the fraction below we obtain: 
\aaa{
\label{22-09}
\frac{\intl_M p(t-\dl,x,t,y) f(y) \la_M(dy)}{\intl_M p(t-\dl,x,t,y) \la_M(dy)}
= \frac{
\intl_{M_t} q(\dl, x_t, y_t)  
e^{-(y_t-x_t,\sg'(t-\te\dl)x)_{\Rnu^m}}f_t(y_t) \la_{M_t}(dy_t)}
{
\intl_{M_t} q(\dl, x_t, y_t) 
e^{-(y_t-x_t,\,\sg'(t-\te\dl)x)_{\Rnu^m}}\la_{M_t}(dy_t)}.
}
Multiplying the numerator and the denominator 
by $\int_{M_t} q(\dl, x_t, y_t)  \la_{M_t}(dy_t)$, and then applying
Corollary \ref{cor1}, we continue \eqref{22-09}:
\mm{
\int_{M} p^M(t-\dl,x,t,y)\, f(y) \la_M(dy) =  f(x) \\
+\dl\, 
\frac{-\frac12\, \lap_{M_t} f_t(x_t)
+(\nab_{M_t} f_t(x_t), \, \sg'(t-\te\dl)x)_{\Rnu^m} 
-f(x)\td R(t,\dl) + \bar R(t,\dl)
}
{
1-\frac{\dl}2\left.\lap_{M_t}e^{-(y_t-x_t,\,\sg'(t-\te\dl)x)_{\Rnu^m}}\right|_{y_t=x_t} 
+ \dl\,\td R(t,\dl)}
}
where $\te\in (0,1)$ is the number satisfying 
$\sg(t)x - \sg(t-\dl)x = \dl\,\sg'(t-\te\dl)x$,
and the functions $\bar R(t,\dl)$ and $\td R(t,\dl)$
are the higher-order terms that appear in the numerator
and the denominator of \eqref{22-09} after applying Corollary \ref{cor1}.
The term $(\nab_{M_t} f_t(x_t), \, \sg'(t-\te\dl)x)_{\Rnu^m}$ appears after computing
\aa{
\left.\nab_{M_t}e^{-(y_t-x_t,\,\sg'(t-\te\dl)x)_{\Rnu^m}}\right|_{y_t=x_t}
= -\Pr\nolimits_{T_{x_t}(M_t)} \sg'(t-\te\dl)x
}
where $\Pr_{T_{x_t}(M_t)}$ denotes the projection onto 
the tangent space ${T_{x_t}(M_t)}$.
Due to the continuity of the map $t \mto \sg'(t)x$, $\sg'(t-\te\dl)x$ converges to 
$\sg'(t)x$ uniformly in $t$ as $\dl\to 0$. Also, as $\dl\to 0$, 
$\frac{\dl}2\left.\lap_{M_t}e^{-(y_t-x_t,\,\sg'(t-\te\dl)x)_{\Rnu^m}}\right|_{y_t=x_t}$
converges to zero uniformly in $t$ by boundedness of the second multiplier.
Therefore, to show that Assumption \ref{limit} is fulfilled we have to prove that
$\bar R(\dl, t)$ and $\td R(t,\dl)$ tend to zero uniformly in $t$ as $\dl \to 0$.
We prove it for the function $\bar R(t,\dl)$.
In the proof of Proposition \ref{pr1} we considered the neighborhoods $U_y=\{z\in M, |z-y|<\eps_y\}$
where the normal coordinates can be introduced. Moreover $\eps_y$ is bounded away from zero by $\eps$ as $y\in M$ varies.
Let $U_{y_t}= \sg(t)U_y$, where $y_t = \sg(t)y$,  and $U_t = \sg(t)U$. Clearly, the exponential map
$\exp: U_t \to U_{y_t}$ is well-defined, and therefore we can introduce normal coordinates in $U_{y_t}$.
Let $\eps_{y_t} = \inf\{|z-y_t|, z\in U_{y_t}\}$.
Due to the continuity of the map $t\mto \sg(t)$, $\eps_{y_t}$ are bounded away from zero, say, by $\eps$,
as $t$ runs over $[S,T]$ and $y$ runs over $M$, i.e. when $y_t$ runs over $\cup_{\tau\in [S,T]} M_\tau$.
This and estimate \eqref{linaaa} imply that
\aa{
 \frac1{(2\pi \tau)^{\frac{d}2}}
 \int_{M_t\dd U_{y_t}} f_t(z) e^{-\frac{|z-y_t|^2}{2\tau}} \la_{M_t}(dz)
 < \tau^{3/2}.
 }
Next, we have to analyze the higher-order term 
in every neighborhood $U_t = \exp^{-1}U_{y_t}$.
We use the estimate of this term obtained in
\cite{sm_vr} (Lemma $2$). All multipliers in the function estimating
the higher-order term
as well as the integral over $\Rnu^m\dd U_t$ are continuous in $t\in [S,T]$.
This proves that $\bar R(t,\dl)$ is bounded by $K\, \dl^{\frac12}$
where $K$ does not depend on $t$.
Now the statement of the theorem follows from Theorem \ref{chernnon}.
\end{proof}
Let us discuss now a probabilistic interpretation of Theorem \ref{thm_app}.
Let $W_\xi$, $\xi\in [\tau-h,\tau]$, be an $\Rnu^m$-valued Brownian motion starting at $0$, 
and let $\W^x_\sg$ be the distribution of the process $x+\sg(\xi)^{-1}W_\xi$.
Further let $U_\eps(M)$ denote the $\eps$-neighborhood of $M$, and
$g:\C([\tau-h,\tau],\Rnu^m) \to \Rnu$ be a $\W^x_\sg$-measurable function.
The right-hand side of the equality
\aa{
\int_{\C([\tau-h,\tau],\Rnu^m)} g(\om) \,\W^x_{\eps,\tau}(d\om) =
\frac{\int_{\C([\tau-h,\tau],\Rnu^m)} \ind_{\{\om:\,\om(\tau)\in U_\eps(M)\}} g(\om) \,\W^x_\sg(d\om)}
{\W^x_\sg \{\om:\om(\tau)\in U_\eps(M)\}}
}
defines a probability distribution $\W^x_{\eps,\tau}$ 
on the same $\sg$-algebra where the distribution $\W^x_{\sg}$
is defined, i.e. 
the $\sg$-algebra generated by all cylindric subsets
of the space of all functions $[\tau-h,\tau]\to \Rnu^m$. Clearly,
$\W^x_{\eps,\tau}$ is supported on $\C([\tau-h,\tau],\Rnu^m)$.
The diffusion associated with $\W^x_{\eps,\tau}$ 
is a time-inhomogeneous Markov process that starts at $x\in M$ at time $\tau-h$, and is
conditioned to come to the neighborhood $U_\eps(M)$ at time $\tau$.
The transition probability function $P_{\tau-h,\tau}(x,\fdot)$ defined via $\W^x_{\eps,\tau}$,
i.e. $P_{\tau-h,\tau}(x,A) = \W^x_{\eps,\tau}(\om: \om(\tau) \in A)$,
possesses the density 
\aa{
p_{\eps}(\tau-h,x,\tau, y) = \frac{\ind_{U_\eps(M)}(y)\,p(\tau-h,x,\tau,y)}{\int_{U_\eps(M)}p(\tau-h,x,\tau,y) dy}.
}
As $\eps$ tends to zero, $p_{\eps}(\tau-h,x,\tau, y)\,dy$ converges weakly relative to the family of
bounded continuous functions to $p^M(\tau-h,x,\tau,y)\,\la_M(dy)$.
The latter function defines a probability distribution on  
the algebra of cylindric subsets of $\C([\tau-h,\tau],\Rnu^m)$. The Markov process associated
to this probability distribution starts at $x\in M$ at time $\tau-h$, and is conditioned
to return to $M$ at time $\tau$.
Consider a partition  
 $\mc P=\{s=t_0 < t_1 < \cdots < t_n=t\}$ of an interval $[s,t]\sub [S,T]$,
and think of a Markov process $X_t^{\mc P}$ that starts at $x\in M$ at time $s$ and is 
conditioned to return to $M$ at all times $t_i \in \mc P$. 
Let $t_{i-1} \lt r < \tau \lt t_i$. 
If $\tau = t_i$ then the transition probability function $P^{\mc P}(r, z,\tau, \fdot)$ of $X_t^{\mc P}$,
considered as a measure, is concentrated on $M$ and  
$p^M(r,z,t_i,y)$ is its density with respect to the measure $\la_M$.
Moreover, the latter holds also if $t_{i-1} < r$.
If $\tau < t_i$ then $P^{\mc P}(r, z,\tau, \fdot)$ is a distribution on the enveloping space $\Rnu^m$.
The conditional probability argument implies that $P^{\mc P}(r, z,\tau, \fdot)$ has the density
with respect to the Lebesgue measure on $\Rnu^m$:
\aaa{
\label{density0}
p^{\mc P}(r,z,\tau,y) = \frac{p(r, z, \tau, y) \int_M p(\tau, y, t_{i},\bar x) \la_M(d\bar x)}
{\int_M p(r, z, t_{i},\bar x) \la_M(d\bar x)}.
}
Now let $t_{i-1} \lt r < t_i < t_{j-1} < \tau < t_j$. In this case 
the density of $P^{\mc P}(z, r,\tau, \fdot)$ with respect to the Lebesgue measure on
$\Rnu^m$ is given by
\mmm{
\label{density-general}
p^{\mc P}(r,z,\tau, y) = 
\int_M p^M(r,z,t_i,x_i) \la_M(dx_i)\int_Mp^M(t_i,x_i,t_{i+1},x_{i+1})\la_M(dx_{i+1})\\
\ldots\int_M p^M(t_{j-2},x_{j-2},t_{j-1},x_{j-1}) \,p^{\mc P}(t_{j-1},x_{j-1},\tau, y)\la_M(dx_{j-1}).
}
\begin{cor}
As the mesh of $\mc P$ tends to zero, 
the finite-dimensional distributions of the process $X_t^{\mc P}$ converge
weakly to the finite-dimensional distributions of the $M$-valued diffusion $X_t$ 
genarated by $A_t$.
\end{cor}
\begin{proof}
We have to prove that for any partition $s< \tau_1 < \cdots < \tau_k < t$ and for any
bounded continuous function $f: \Rnu^k \to \Rnu$,
\aa{
\E\bigl[f(X^{\mc P}_{\tau_1}, \ldots, X^{\mc P}_{\tau_k})\bigr]
\to 
\E\bigl[f(X_{\tau_1}, \ldots, X_{\tau_k})\bigr]
\quad \text{as} \; |\mc P| \to 0. 
}
First we consider only those partitions $\mc P$ that contain all the points $\tau_i$, $1\lt i \lt k$. 
Pick up two subsequent points $\tau_i$ and $\tau_{i+1}$. Let $t_l\in \mc P$ and $t_m\in \mc P$ be 
such that $t_l = \tau_i$ and $t_m = \tau_{i+1}$.  Then,
\mmm{
\label{density-particular}
p^{\mc P}(\tau_i,z,\tau_{i+1}, y) = 
\int_M p^M(\tau_i,z,t_{l+1},x_{l+1}) \la_M(dx_{l+1})\\
\ldots\int_M p^M(t_{m-2},x_{m-2},t_{m-1},x_{m-1}) \,p^{M}(t_{m-1},x_{m-1},\tau_{i+1}, y)\la_M(dx_{m-1}).
}
Now let $f\in \C(M^k)$, and let $P(s,x,t,A)$ be the transition
probability function of the process $X_t$ on $M$ generated by $A_t$. We have:
\mm{
\int_M p^{\mc P}(s,x,\tau_1,x_1)\,\la_M(dx_1) 
\ldots \int_M p^{\mc P}(\tau_{k-1},x_{k-1},\tau_k,x_k)f(x_1,\ldots, x_k)\,\la_M(dx_k)
\\
- \int_M P(s, x, \tau_1, dx_1) 
\ldots \x \int_M P(\tau_{k-1}, x_{k-1}, \tau_k, dx_k) f(x_1,\ldots, x_k)=
}
\mmm{
\label{technique0}
= \sum_{i=1}^{k-1}  \int_M p^{\mc P}(s,x,\tau_1,x_1)\,\la_M(dx_1) 
\ldots \int_M p^{\mc P}(\tau_{i-1},x_{i-1},\tau_i,x_i)\,\la_M(dx_i)\\
\Bigl[\int_M p^{\mc P}(\tau_i, x_i, \tau_{i+1}, x_{i+1})\la_M(dx_{i+1})
 - P(\tau_i,x_i, \tau_{i+1}, dx_{i+1})
\Bigr]\\
\int_M P(\tau_{i+1}, x_{i+1}, \tau_{i+2}, dx_{i+2})
\ldots \int_M P(\tau_{k-1},x_{k-1},\tau_k, dx_k)\,
f(x_1,\ldots, x_k)\,\la_M(dx_k).}
Each term of this sum converges to zero as the mesh $|\mc P|$ goes to zero.
Indeed, for every $i$, $1\lt i < k$, for every function $g\in \C(M^{i+1})$,
the difference
\mmm{
\label{chernoff-difference}
\int_M \bigl(p^{\mc P}(\tau_i, \td x_i, \tau_{i+1}, x_{i+1})\, 
\la_M(dx_{i+1}) 
 - \int_M P(\tau_i, \td x_i, \tau_{i+1}, dx_{i+1})\bigr)
g(x_1, \ldots, x_{i+1})
}
converges to zero.
This follows from Theorem \ref{thm_app}. Indeed, 
the second term in \eqref{chernoff-difference} 
is the backward propagator with the left generator $A_t$, 
and the first term is the operator product \eqref{Qprod}. 
The convergence in \eqref{chernoff-difference} holds in $\C(M^{i+1})$
by the argument of Theorem \ref{thm_app}. The latter argument has to be 
applied to operators  $\C(M^{i+1})\to \C(M^{i+1})$ and with respect to the norm of
$C(M^{i+1})$ instead of $\C(M)$, as in Theorem \ref{thm_app}, 
which, however, leaves the proof of Theorem \ref{thm_app} without changes.

Now let us assume that an infinite number of partitions $\mc P$
with the meshes decreasing to zero do not include some of the points
$\tau_i$.
Then, instead of formula \eqref{density-particular} for $p^{\mc P}$ 
we have to use formula \eqref{density-general}.
We would like to reduce this case to the previous one, i.e. when 
all the points $\tau_i$ are always among the partition points of $\mc P$.
For this purpose we have to analyze the expression:
\aaa{
\label{27}
\int_{\Rnu^m} p^{\mc P}(t_i, x_i, \tau, y)\, dy \int_M p^{M}(\tau, y, x_{i+1}, t_{i+1})
f(z, y, x_{i+1})\, \la_M(dx_{i+1}),
}
where the variable $x_{i+1}$ comes from the subsequent integrals as the result of their
replacement with the backward propagator. The variables $z$ and $y$ come from the original 
integrand function,
$\tau$ is one of the points $\tau_i$, $1\lt i \lt k$, and
the partition points $t_i$ and $t_{i+1}$ are choosen such that
$t_i < \tau < t_{i+1}$.
We would like to show that as $t_i, t_{i+1} \to \tau$, the difference between 
\eqref{27} and $\int_M p^M(t_i,x_i,t_{i+1}, x_{i+1})\,f(z,x_i,x_{i+1})\, \la_M(dx_{i+1})$
converges to zero. By the Banach-Steinhaus
theorem, it suffices to prove this when $f$ is 
continuously differentiable with respect to $y$. 
Applying formula \eqref{density0} we observe that
\eqref{27} equals to
\aaa{
\label{fraction0}
\frac{\int_{\Rnu^m} dy \, p(t_i, x_i, \tau, y)
\int_M p(\tau, y, t_{i+1}, x_{i+1})\,  f(z, y, x_{i+1})\,\la_M(dx_{i+1})}
{\int_M p(t_i, x_i, t_{i+1}, \bar x_{i+1})\, \la_M(d\bar x_{i+1})}.
}
Applying formula \eqref{int.dens}, 
for the numerator of \eqref{fraction0} we obtain:
\aaa{
\label{analyze}
\int_{\Rnu^m} \hspace{-2mm} dy_{\tau} \, q(\tau-t_i, x_{t_i}, y_{\tau})
\int_{M_{t_{i+1}}} \hspace{-4mm} q(t_{i+1} - \tau, y_{\tau}, x_{t_{i+1}})\,
f_{\tau,t_{i+1}}(z,y_\tau, x_{t_{i+1}}) 
\la_{M_{t_{i+1}}}\!(dx_{t_{i+1}})
}
where $f_{\tau,t_{i+1}}(z,\fdot,\fdot) = 
f\bigl(z, \sg(\tau)^{-1}(\fdot), \sg(t_{i+1})^{-1}(\fdot)\bigr)$,
$x_{t_i} = \sg(t_i)x_i$, $y_\tau = \sg(\tau)y$,
$x_{t_{i+1}} = \sg(t_{i+1})x_{i+1}$.
Also, in \eqref{analyze} we omitted the multiplier ${\det{\sg(t_{i+1})}}$
which will be taken into consideration later again.
Application of Taylor's formula to $f_{\tau,t_{i+1}}(z,\fdot,x_{t_{i+1}})$
at point $x_{t_i}$ gives:
\aaa{
\label{taylor}
f_{\tau,t_{i+1}}(z,y_\tau,x_{t_{i+1}}) = f_{\tau,t_{i+1}}(z,x_{t_i},x_{t_{i+1}})
+ \pl_2 f_{\tau,t_{i+1}}(z,p(x_{t_i},y_\tau),x_{t_{i+1}})(y_\tau - x_{t_i})
}
where $p(x_{t_i},y_\tau)$ is a point on the segment $[x_{t_i},y_\tau]$
and $\pl_2$ means partial differentiation with respect to the second argument.
If we substitute as an integrand the first summand of \eqref{taylor} into \eqref{analyze}, 
we obtain:
\aa{
\int_{M_{t_{i+1}}} q(t_{i+1}-t_i, x_{t_i}, x_{t_{i+1}}) f_{\tau,t_{i+1}}(z,x_{t_i},x_{t_{i+1}})
\la_{M_{t_{i+1}}}\!(dx_{t_{i+1}}).
}
Further, this substitution brings \eqref{fraction0} to
\aa{
\int_M p^M(t_i,x_i,t_{i+1},x_{i+1}) \, f(z,\sg(\tau)^{-1}\sg(t_i)x_i, x_{i+1})
\, \la_M(dx_{i+1}).
}
The latter converges to 
$\int_M p^M(t_i,x_i,t_{i+1},x_{i+1}) \, f(z,x_i, x_{i+1})\, \la_M(dx_{i+1})$.
Thus, we have to prove that 
\aaa{
\label{1-10}
&\frac{\det\sg(t_{i+1})}{\int_M p(t_i,x_i, t_{i+1}, x_{i+1})\la_M(dx_{i+1})}
\int_{\Rnu^m}  dy_{\tau} \, q(\tau-t_i, x_{t_i}, y_{\tau})  \\ 
&\x
\int_{M_{t_{i+1}}} \hspace{-2mm} q(t_{i+1} - \tau, y_{\tau}, x_{t_{i+1}})\,
\pl_2 f_{\tau,t_{i+1}}(z,p(x_{t_i},y_\tau),x_{t_{i+1}})(y_\tau - x_{t_i})\,
\la_{M_{t_{i+1}}}\!(dx_{t_{i+1}})\notag
} 
converges to zero as $t_i, t_{i+1} \to \tau$.
We change the order of integration in \eqref{1-10} and split the integral
with respect to $y_\tau$, taken over $\Rnu^m$, into
two: over the set $\{y_\tau: \, |y_\tau - x_{t_i}| < (\tau-t_i)^{\frac13}\}$
and over its complement $\{y_\tau: \, |y_\tau - x_{t_i}|\gt (\tau-t_i)^{\frac13}\}$.
Estimation of the first term gives:
\aaa{
&\Bigl|\int_{\{|y_\tau - x_{t_i}| < (\tau-t_i)^{\frac13}\}} \hspace{-1mm} dy_{\tau} \, q(\tau-t_i, x_{t_i}, y_{\tau})
\, q(t_{i+1} - \tau, y_{\tau}, x_{t_{i+1}}) \notag\\
&\phantom{\Bigl|\int_{\{|y_\tau - x_{t_i}| < (\tau-t_i)^{\frac13}\}}}
\x\pl_2 f_{\tau,t_{i+1}}(z,p(x_{t_i},y_\tau),x_{t_{i+1}})(y_\tau - x_{t_i})\Bigr|
\lt \label{4-10}\\
&\sup_{x_i, x_{i+1}\in M, y\in U_\eps(M),z\in K}|\pl_2 f_{\tau,t_{i+1}}(z,p(x_{t_i},y_\tau),x_{t_{i+1}})|
\,q(t_{i+1} - t_i,x_{t_i}, x_{t_{i+1}})\,(\tau-t_i)^\frac13\notag
}
where $K$ is a compact,
since without loss of generality we can consider that the totality of variables $z$
belongs to a compact $K$. Indeed, the integrand 
$\pl_2 f_{\tau,t_{i+1}}(z,p(x_{t_i},y_\tau),x_{t_{i+1}})(y_\tau - x_{t_i})$ 
is bounded by \eqref{taylor}. The preceeding integration w.r.t. every variable 
from the totality $z$ is either taken over the manifold $M$ or similar to 
the integration w.r.t. $y$ in \eqref{27}. In the latter case we can replace the integrals
over $\Rnu^m$ with integrals over compact neighborhoods of $M$ because 
the integrals over  the complements of these neighborhoods tend to zero as 
the mesh $|\mc P|$ goes to zero.  
Further, since $y_\tau$ is always in the $(\tau-t_i)^\frac13$-neighborhood of $M_{t_i}$,
then $y$ is always in some $\eps$-neighborhood of $M$. 
Therefore, the supremum in the last line of \eqref{4-10} is finite.
Hence, the summand in \eqref{1-10} which corresponds to the integration over
$\{y_\tau: \, |y_\tau - x_{t_i}| < (\tau-t_i)^{\frac13}\}$
is bounded by
\aa{
\sup|\pl_2f_{\tau,t_{i+1}}|\, (\tau-t_i)^\frac13 \int_M p^M(t_i, x_i, t_{i+1}, x_{i+1})
\la_M(dx_{i+1})
}
which tends to zero as $t_i, t_{i+1}\to \tau$.
The other summand, which corresponds to the integration over
$\{y_\tau: \, |y_\tau - x_{t_i}| \gt (\tau-t_i)^{\frac13}\}$,
is bounded by
\mm{
\frac1{\bigl(2\pi(\tau-t_i)\bigr)^{\frac{m}2}}e^{-\frac1{2(\tau-t_i)^\frac13}}
\int_{M_{t_{i+1}}}\la_{M_{t_{i+1}}}(dx_{t_{i+1}})\\
\x\int_{\Rnu^m} dy_\tau \bigl|\pl_2 
f_{\tau,t_{i+1}}(z,p(x_{t_i},y_\tau),x_{t_{i+1}})(y_\tau - x_{t_i})\bigr|
\,q(t_{i+1}-\tau, y_\tau,x_{t_{i+1}})
}
which also tends to zero as $t_i, t_{i+1} \to \tau$ since
the product of the partial derivative $\pl_2 f_{\tau,t_{i+1}}$ and $(y_\tau-x_{t_i})$ is
bounded by formula \eqref{taylor}.  The multiplier in front of the integrals in \eqref{1-10}
also converges to zero as $|\mc P| \to 0$, and therefore
the convergence of \eqref{1-10} to zero as $t_i,t_{i+1}\to \tau$ is proved. 
\end{proof}
\subsection*{Acknowledgements} 
The author would like to thank the referee for meaningful comments.
This work was supported by the Portuguese Foundation for 
Science and Technology through  
the Centro de Matem\'atica da Universidade do Porto.

\end{document}